\documentclass{article}
\usepackage{amssymb}
\newtheorem{teo}{Theorem}
\newtheorem{lem}{Lemma}

\newcommand{\ri}{{\rm ri\,}}

\newcommand{\cl}{{\rm cl\,}}
\newcommand{\CL}{{\rm CL\,}}
\newcommand{\cone}{{\rm cone\,}}
\newcommand{\intern}{{\rm int\,}}
\newcommand{\conv}{{\rm conv\,}}
\newcommand{\dom}{{\rm dom\,}}

\begin{document}
\title{Constructive no-arbitrage criterion under transaction costs
in the case of finite discrete
time\footnote{Russian version submitted to Theory of Probability and
its Applications: November 2005}}
\author{Dmitry B. Rokhlin}
\date{}
\maketitle
\begin{abstract}
We obtain a constructive criterion for robust no-arbitrage in
discrete-time market models with transaction costs.
This criterion is expressed in terms of the supports of
the regular conditional upper distributions
of the solvency cones. We also consider the model with a
bank account. A method for construction of arbitrage strategies
is proposed.
\end{abstract}
{\bf Introduction.} According to the classical Dalang-Morton-Willinger
theorem the perfect market is arbitrage-free iff there exists
an equivalent martingale measure for the discounted price process
of traded assets \cite{DMW90}. In the paper
\cite{JS98} there was explicitly formulated one more
equivalent condition, expressed in terms of the supports
of the regular conditional distributions of price
increments (see condition (g) of Theorem 3 in the cited paper).
It is natural to call this condition {\it constructive} since it
gives a computationally feasible method for the verification
of the absence of arbitrage. An analogous condition under
convex portfolio constraints was introduced in \cite{R04}.

The aim of the present paper is to obtain a similar
constructive criteria for no-arbitrage in the framework of the market
models with transaction costs introduced in the papers \cite{K99},
\cite{KS01} and further developed in \cite{KRS02}, \cite{S04}, \cite{KRS03},
\cite{G05}. For these models {\it strictly consistent price processes}
are considered instead of equivalent martingale measures.
This name is assigned to martingales, taking their values in
the relative interior of the random cones $K_t^*$, conjugate to
the solvency cones \cite{KS01}, \cite{KRS02}, \cite{S04}.
In the paper \cite{S04} there was introduced the condition
of the {\it robust no-arbitrage} (NA$^r$) and was proved that the
existence of a strictly consistent price process is equivalent to the
fulfillment of this condition.

Condition (iii) of the paper \cite{KRS03} gives a purely geometrical
characterization of the robust no-arbitrage.
In Section 1 of the present paper, using the notion of the
{\it support of the regular conditional upper distribution}
of a measurable set-valued map \cite{R05}, we give the dual
"constructive" description of the mentioned condition [8]
(see Theorem 1 below). In Section 2, where the market model
with transaction costs is considered, we show that
this result leads to a desirable no-arbitrage criterion
(Theorem 2, condition (c)).

We mention that under the assumption on the existence of
interior points of $K_t^*$ (known as the {\it efficient friction
condition} \cite{KRS02}), and also for the case of a finite probability
space, this criterion, in fact, was obtained in \cite{R05}.

In Section 3 we consider an important partial case: the market
model with a bank account. Here the known no-arbitrage criterion
is formulated in a more traditional way and states the existence
of a price process, lying between the bid and ask prices,
and a correspondent equivalent probability measure (cf.
\cite{JK95}). The constructive no-arbitrage criterion is expressed
in terms of the Cartesian products of the bid-ask price intervals
(Theorem 3).

The violation of NA$^r$ condition means the
existence of an arbitrage strategy for any bid-ask process
with lower transaction costs (the precise formulation is given in
\cite{S04} and in Section 2 of the present paper).
As is well known, on a perfect market arbitrage strategies
involve only two changes of the portfolio in neighbour time moments.
They are easily constructed by the separation arguments \cite{JS98}.
However, under transaction costs such strategies can have a more
complex structure. The reasoning given below, enable to extract
some information about them. Particularly, it is possible to
determine the time moment for taking the first non-trivial portfolio
and to specify how to do it. An illustrative example
is given in concluding Section 4.

{\bf 1. Key result.}
\setcounter{section}{1}
\setcounter{equation}{0}
Denote by $\cl A$, $\intern A$, $\ri A$, $\conv A$, $\cone A$
the closure, the interior, the relative interior,
the convex hull, and the conic hull of a subset $A$ of a finite dimensional
space \cite{R73}. Put also
$$ A\pm B=\{x\pm y: x\in A,\ y\in B\},\ \ \
   \lambda A=\{\lambda x:x\in A\}.$$
If $A$ is a cone (i.e. $\lambda A\subset A$, $\lambda\ge 0$),
then $A^*$ is the conjugate cone:
$A^*=\{y\in\mathbf R^d:\langle x,y\rangle\ge 0,\ x\in A\}$.
Here $\langle x,y\rangle$ is the usual scalar product of the vectors
$x$ and $y$. Let
$A^\perp=\{y:\langle x,y\rangle=0,\ x\in A\}$ and
$$ H^x=(\cone x)^\perp,\ \ \ H^x_\ge=(\cone x)^*$$
for any $x\in\mathbf R^d$.

Consider a probability space $(\Omega,\mathcal F, \mathbf P)$ and
a $\sigma$-algebra $\mathcal H\subset\mathcal F$.
In the sequel we assume all $\sigma$-algebras to be
complete with respect to the measure $\mathbf P$.

A set-valued map $F$, assigning some set $F(\omega)\subset \mathbf R^d$
to each $\omega\in\Omega$, is called $\mathcal H$-measurable if
$\{\omega:F(\omega)\cap V\neq\emptyset\}
\in\mathcal H$ for any open set $V\subset \mathbf R^d$.
The {\it effective set} of $F$ is defined as follows:
$\dom F=\{\omega:F(\omega)\neq\emptyset\}.$
The function $\xi:\Omega\mapsto \mathbf R^d$ is called
a {\it selector} of $F$ if $\xi(\omega)\in F(\omega)$
for all $\omega\in\dom F$.

A mapping $F$ is called closed (resp. convex, conic) if
the sets $F(\omega)$ are closed (resp. are convex, are cones). The
$\mathcal H$-measurability of the closed set-valued map
$F$ is equivalent to the existence of a countable family
$\{\xi_i\}_{i=1}^\infty$ of $\mathcal H$-measurable selectors of $F$
such that the sets $\{\xi_i(\omega)\}_{i=1}^\infty$
are dense in $F(\omega)$ for all $\omega\in\dom F$ \cite[p.331]{L85}.
The mentioned family is called a {\it Castaing representation} for $F$.

Let $\CL=\CL(\mathbf R^d)$ be the family of nonempty closed subsets
of $\mathbf R^d$. Denote by $\mathcal E(\CL)$ the {\it
Effros $\sigma$-algebra} \cite{B91}, \cite{S98}, generated by
the sets
$$A_V=\{D\in\CL: D\cap V\neq\emptyset\},$$
where $V$ is an open subset of $\mathbf R^d$.
In the paper \cite{R05} there was mentioned that for an
$\mathcal F$-measurable set-valued map $F$ with nonempty closed
values there exists a function
$$ \mathbf P^*:\Omega\times\mathcal E(\CL)\mapsto\mathbf R$$
with the following properties:
\begin{itemize}
\item[\rm (i)] for any $\omega$ the function $C\mapsto\mathbf P^*(\omega,C)$
is a probability measure on $\mathcal E(\CL)$;
\item[\rm (ii)] for any $C\in\mathcal E(\CL)$ the function
$\omega\mapsto\mathbf P^*(\omega,C)$ a.s. coincides with
$\mathbf P(\{F\in C\}|\mathcal H)(\omega)$.
\end{itemize}

This result allows to introduce the {\it regular conditional
upper distribution} of a set-valued map and the support of this
distribution \cite{R05}. The function $\mu_F(\omega,V)=P^*(\omega,A_V)$,
where $V$ is an open subset of $\mathbf R^d$, is called
the  regular conditional upper distribution of $F$ with respect
to $\mathcal H$. The set
$$ \mathcal K(F,\mathcal H;\omega)=
\{x\in\mathbf R^d:\mu_F(\omega,B_\varepsilon(x))>0,
 \ \forall \varepsilon>0\}$$
is called the support of $\mu_F(\omega,\cdot)$

The set-valued map $\omega\mapsto K(F,\mathcal H;\omega)$
takes nonempty values and is $\mathcal H$-measurable
and closed \cite[Proposition\,4(a)]{R05}.

If the set $\Omega$ is finite, the probability measure is non-degenerate
(i.e. $\mathbf P(\omega)>0$, $\omega\in\Omega$), and the
algebra $\mathcal H\subset \mathcal F$ is generated by the
partition $\{D_1,\dots,D_l\}$, then, as it mentioned in \cite{R05}, we have
\begin{equation} \label{1.1}
\mathcal K(F,\mathcal H;\omega)=\bigcup_{\omega'\in D_i} F(\omega'),
\ \ \omega\in D_i,\ \ i=1,\dots, l.
\end{equation}

If the mapping $F$ is single-valued, then $\mu_F$ and
$\mathcal K(F,\mathcal H)$ are respectively the regular conditional
distribution of $F$ with respect to $\mathcal H$ and its support.

The mapping $\mathcal K(F,\mathcal H)$ can be defined also if
$F$ have the empty values on a set of positive measure.
For that consider the closed extension $F_*$ of the mapping $F$
to $\Omega$ (for example,
$F_*=FI_{\{F\neq\emptyset\}}+\{0\}I_{\{F=\emptyset\}}$). Up
to a zero measure set, the mappings $\mathcal K(F_*,\mathcal H)$
coincide for all such extensions \cite[Proposition\,4(b)]{R05}.
We regard $\mathcal K(F,\mathcal H)$ as a representative of
the correspondent equivalence class.

If $F(\omega)=\emptyset$ on a set of positive measure, then we put
$\mathcal K(F,\mathcal H)=\emptyset$ for all $\omega$.

Let $\{\xi_i\}_{i=1}^\infty$ be a Castaing representation for
an $\mathcal F$-measurable map $F:\Omega\mapsto\CL$.
In the sequel we use the following result (\cite[Lemma 1]{R05}):
\begin{equation} \label{1.2}
 \mathcal K(F,\mathcal H)=\cl\left(\bigcup_{i=1}^\infty
 \mathcal K(\xi_i,\mathcal H)\right) \ \ {\rm a.s.}
\end{equation}
Certainly, this relation still holds true in the case where $F$
have the empty values on a negligible set.

Note that if a closed mapping $F$ is conic, then
the mapping $\mathcal K(F,\mathcal H)$ is conic too.
Indeed, if $\{\xi_i\}_{i=1}^\infty$ is a Castaing representation for
$F$, then the family
$\{\lambda \xi_i\}_{i=1}^\infty$, $\lambda>0$
possesses the same property. Using formula
(\ref{1.2}), we get
$$ \mathcal K(F,\mathcal H)={\rm cl}\left(\bigcup_{i=1}^\infty
   \mathcal K(\lambda\xi_i,\mathcal H)\right)=
\lambda\,{\rm cl}\left(\bigcup_{i=1}^\infty
\mathcal K(\xi_i,\mathcal H)\right)
= \lambda \mathcal K(F,\mathcal H).$$

Now we turn to the formulation of the key result of this paper.
Suppose we have a sequence $\mathcal F_t$-measurable
mappings $G_t:\Omega\mapsto \CL$, $t=0,\dots,T$, whose values
are the closed convex cones. We introduce the sequence of the set-valued
mappings $W_t$ by the recursion:
$$ W_T=G_T^*,$$
$$ W_t=\cl(\ri G_t^*\cap \ri Y_t),\ \
   Y_t=\cl(\conv\mathcal K(W_{t+1},\mathcal F_t)),\ \ \
   0\le t\le T-1.$$

Using Proposition 3 of the paper \cite{R05} and Lemma 1 of
\cite{R04}, it is not difficult to check that
the $\mathcal F_{t+1}$-measurability of the set-valued map
$W_{t+1}$ imply the $\mathcal F_t$-meaurability of
$W_t$. Consequently, the sequence $W_t$ is well-defined.
If $W_{n+1}=\emptyset$ on a set of positive measure, then
$W_t=\emptyset$ for $t\le n$.

Denote by $L^0(F,\mathcal H)$ the set of all $\mathcal H$-measurable
functions $\xi$, satisfying the condition $\xi(\omega)\in F(\omega)$
for almost all $\omega$.
\begin{teo} The following conditions are equivalent:
\begin{itemize}
\item[\rm (a)] if $\sum_{t=0}^T x_t=0$, where
$x_t\in L^0(-G_t,\mathcal F_t)$, then $x_t=L^0(G_t\cap(-G_t),\mathcal F_t)$,
$t=0,\dots,T$;
\item[\rm (b)] $W_t\neq\emptyset$ a.s., $t=0,\dots,T-1$.
\end{itemize}
\end{teo}

This theorem gives the dual description of condition (a),
which is borrowed from the papers \cite{KRS02}, \cite{KRS03}.
Note, that in the case of a one-point probability space we have
$Y_t=W_{t+1}$ and $W_t=\cl (\ri G_t^*\cap \ri W_{t+1})$.
Here condition (b) shapes to the form
$\cap_{t=0}^T \ri G_t^*\neq\emptyset$
(taking into account Theorem 6.5 of \cite{R73}).

The proof of Theorem 1 is forestalled by a number
of auxiliary statements.
\begin{lem} For a conic $\mathcal F$-measurable
mapping $F:\Omega\mapsto \CL$
the following equalities hold true
\begin{equation} \label{1.3}
 L^0(F^*,\mathcal H)=
 L^0\left(\mathcal K(F,\mathcal H)^*,\mathcal H\right),
   \ \ L^0(F^\perp,\mathcal H)=
 L^0\left(\mathcal K(F,\mathcal H)^\perp,\mathcal H\right).
\end{equation}
\end{lem}
{\it Proof.} Let $\{\xi_i\}_{i=1}^\infty$ be a Castaing representation
for $F$. Evidently, the condition $x\in L^0(F^*,\mathcal H)$ means
that $\mathbf P (\xi_i\in H_\ge^x)=1$ for all $i\ge 1$.
By Lemma 3 of \cite{R04} we have
$$ \mathbf P(\xi_i\in H_\ge^x)=\mathbf E\mathbf P(\xi_i\in H_\ge^x
|\mathcal H)=\int_\Omega\mu_{\xi_i}(\omega,H_\ge^{x(\omega)})\,
d\mathbf P(\omega)=1,$$
where $\mu_{\xi_i}$ is the regular conditional distribution of $\xi_i$
with respect to $\mathcal H$.
This equality is equivalent to the relation
$\mu_{\xi_i}(\omega,H_\ge^{x(\omega)})=1$ a.s., which,
in turn, means that $\mathcal K(\xi_i,\mathcal H)\subset
H_\ge^{x(\omega)}$ a.s.

Thus, the condition $x\in L^0(F^*,\mathcal H)$
is equivalent to the following set of conditions:
$x\in L^0((\cone \mathcal K(\xi_i,\mathcal H))^*,
\mathcal H)$, $i\ge 1$. The first equality (\ref{1.3})
is now implied by formula (\ref{1.2}):
$$\mathcal K(F,\mathcal H)^*=\bigcap_{i=1}^\infty
 (\cone\mathcal K(\xi_i,\mathcal H))^* \ \ {\rm a.s.}$$
The second equality (\ref{1.3}) is verified in a quite similar
manner.
\begin{lem} \label{L2}
Assume that the $\mathcal H$-measurable maps
$F_1, F_2, F:\Omega\mapsto\CL$
satisfy the condition $F_1+F_2=F$ a.s. Then
$$L^0(F_1,\mathcal H)+L^0(F_2,\mathcal H)=L^0(F,\mathcal H).$$
\end{lem}
{\it Proof.} The inclusion $L^0(F_1,\mathcal H)+
L^0(F_2,\mathcal H)\subset L^0(F,\mathcal H)$ is obvious.
To prove the converse inclusion consider an arbitrary element
$z\in L^0(F,\mathcal H)$ and introduce the set-valued map
$$\omega\mapsto\Psi(\omega)=(F_1(\omega)\times F_2(\omega))
\cap\{(x_1,x_2):x_1+x_2=z(\omega)\}.$$
The map $\Psi$ takes nonempty closed values a.s.
and is $\mathcal H$-measurable as the intersection of closed
$\mathcal H$-measurable maps (its enough to consider the graphs
of the correspondent maps and to make use of the completeness of
the $\sigma$-algebra $\mathcal H$: see \cite{H75} or
\cite[Proposition 3]{R04}).
Any element $\xi=(\xi_1,\xi_2)\in L^0(\Psi,\mathcal H)$
possesses the following properties:  $\xi_1+\xi_2=z$,
$\xi_i\in F_i$ a.s., $i=1,2$. This completes the proof.

\begin{lem} Let $W_t\neq\emptyset$ a.s., $t=n,\dots,T-1$.
Then
$$ L^0(W_n^*,\mathcal F_n)\subset\sum_{t=n}^T
   L^0(G_t,\mathcal F_t).$$
\end{lem}
{\it Proof.} Since the relative interiors of the closed convex cones
$G_n^*$ and $Y_n$ have a common point a.s., it follows that
$W_n=G_n^*\cap Y_n$ (\cite[Theorem 6.5]{R73}) and
$$ W_n^*=G_n+Y_n^*\ \ {\rm a.s.}$$
(\cite[Theorem 16.4.2]{R73}). Using Lemmas 1 and 2, we get
$$ L^0(Y_n^*,\mathcal F_n)=L^0(\mathcal K(W_{n+1},\mathcal F_n)^*,
\mathcal F_n)=L^0(W_{n+1}^*,\mathcal F_n);$$
$$ L^0(W_n^*,\mathcal F_n)=L^0(G_n,\mathcal F_n)+
   L^0(W_{n+1}^*,\mathcal F_n)\subset L^0(G_n,\mathcal F_n)+
   L^0(W_{n+1}^*,\mathcal F_{n+1}). $$
These relations imply the statement of Lemma 3 in an obvious way.

\begin{lem} \label{L4}
Suppose $\sum_{t=0}^T x_t=0$, $x_t\in L^0(-G_t,\mathcal F_t)$,
$t=0,\dots,T$ and $W_t\neq\emptyset$ a.s., $t=n,\dots,T-1$, where $n\ge 1$.
Then
$$ y_n=-\sum_{t=n}^T x_t\in L^0(W_n^*,\mathcal F_{n-1})
   \cap L^0(Y_{n-1}^*,\mathcal F_{n-1}).$$
\end{lem}
{\it Proof.} The function $y_T=-x_T=\sum_{t=0}^{T-1} x_t$
is $\mathcal F_{T-1}$-measurable and satisfy the condition
$y_T\in G_T=W_T^*$ a.s. Using Lemma 1, we obtain
\begin{equation} \label{1.4}
 y_T\in L^0(\mathcal K(W_T,\mathcal F_{T-1})^*,\mathcal F_{T-1}) =
 L^0(Y_{T-1}^*,\mathcal F_{T-1}).
\end{equation}

Thus, the statement of Lemma holds true for $n=T$.
Under the assumption of the validity of this statement for $n=m+1$
let us prove it for $n=m$. As long as $W_m\subset G_m^*$,
it follows that $-x_m\in G_m\subset W_m^*$ a.s.
Furthermore, by the above assumption, $y_{m+1}\in Y_m^*
\subset W_m^*$ a.s. Consequently,
$y_m\in L^0(W_m^*,\mathcal F_{m-1})$, since the function
$y_m=-x_m+y_{m+1}=\sum_{t=0}^{m-1} x_t$ is
$\mathcal F_{m-1}$-measurable.

 By Lemma 1 this yields that
\begin{equation} \label{1.5}
 y_m\in L^0(\mathcal K(W_m,\mathcal F_{m-1})^*,\mathcal F_{m-1})=
   L^0(Y_{m-1}^*,\mathcal F_{m-1}).
\end{equation}

\begin{lem} \label{L5}
Suppose $\sum_{t=0}^T x_t=0$,
$x_t\in L^0(G_t\cap (-G_t),\mathcal F_t)$, $t=0,\dots,T$ and
$W_t\neq\emptyset$ a.s., $t=n,\dots,T-1$, where $n\ge 1$. Then
$$ y_n=-\sum_{t=n}^T x_t\in L^0(W_n^\perp,\mathcal F_{n-1})
   \cap L^0(Y_{n-1}^\perp,\mathcal F_{n-1}). $$
\end{lem}
{\it Proof}. The function $y_T=-x_T=\sum_{t=0}^{T-1} x_t$
is $\mathcal F_{T-1}$-measurable and satisfy the condition
$y_T\in G_T\cap (-G_T)=W_T^\perp$ a.s.
The inclusion $y_T\in L^0(Y_{T-1}^\perp,\mathcal F_{T-1})$ is verified
by Lemma 1 similar to (\ref{1.4}).

Assume that the assertion o Lemma is true for $n=m+1$. Let us prove it
for $n=m$. Since $W_m\subset G_m^*$, we have $-x_m\in G_m\cap (-G_m)\subset
W_m^\perp$ a.s. Furthermore, according to the above assumption,
$y_{m+1}\in Y_m^\perp\subset W_m^\perp$ a.s. This yields that
$y_m\in L^0(W_m^\perp,\mathcal F_{m-1})$, since
the function $y_m=-x_m+y_{m+1}=\sum_{t=0}^{m-1} x_t$ is
$\mathcal F_{m-1}$-measurable.
The inclusion $y_m\in L^0(Y_{m-1}^\perp,\mathcal F_{m-1})$
is now verified similar to (\ref{1.5}).

\begin{lem} Let $A$, $B$, $C$ be convex subsets of
$\mathbf R^d$ such that
$$ C=\cl(\ri A\cap\ri B)\neq\emptyset.$$
Then the sets $A$ and $C$ cannot be separated.
\end{lem}
{\it Proof.} Note that
$$ \ri C\subset\ri(\cl(\ri A))=\ri A.$$
Consequently, $\ri C\cap\ri A=\ri C\neq\emptyset$ and the assertion
of Lemma is implied by the separation criterion for
the convex sets \cite{R73}.

\begin{lem} \label{L7}
Suppose $\sum_{t=0}^T x_t=0$, $x_t\in L^0(-G_t,\mathcal F_t)$
and $W_t\neq\emptyset$ a.s., $t=0,\dots,T-1$. Then
$$ y_n=-\sum_{t=n}^T x_t\in L^0(W_n^\perp,\mathcal F_{n-1}),\ \
   n=0,\dots,T, $$
where $\mathcal F_{-1}=\mathcal F_0$.
\end{lem}
{\it Proof.} Assume that
$ y_n=-\sum_{t=n}^T x_t\not\in L^0(W_n^\perp,\mathcal F_{n-1}), $
for some $n$. Then by Lemma 1 we get
$$ y_n\not\in L^0(\mathcal K(W_n,\mathcal F_{n-1})^\perp,\mathcal F_{n-1})
  =L^0(Y_{n-1}^\perp,\mathcal F_{n-1}). $$
On the other hand,
$$ y_n\in L^0(Y_{n-1}^*,\mathcal F_{n-1})\subset
          L^0(W_{n-1}^*,\mathcal F_{n-1}), $$
by Lemma 4. Therefore, on a set $A_{n-1}\in\mathcal F_{n-1}$
of positive measure we have $y_n\in Y_{n-1}^*\backslash Y_{n-1}^\perp$.

If $y_n\in W_{n-1}^\perp$ on a set $A'_{n-1}\subset A_{n-1}$
of positive measure, then $y_n$ properly separates the sets $W_{n-1}$ and
$Y_{n-1}$ for $\omega\in A'_{n-1}$. But this is impossible by Lemma 6
and the definition of $W_{n-1}$, since $W_{n-1}\neq\emptyset$ a.s.

Thus $y_n\in W_{n-1}^*\backslash W_{n-1}^\perp$ a.s. on $A_{n-1}$.
However, $-x_{n-1}\in G_{n-1}\subset W_{n-1}^*$ a.s. for any $n\ge 1$.
Consequently,
$$y_{n-1}=-x_{n-1}+y_n\not\in L^0(W_{n-1}^\perp,\mathcal F_{n-2}).$$
It follows by induction that $y_0=\sum_{t=0}^T
x_t\not\in L^0(W_0^\perp,\mathcal F_0)$ and we get the contradiction since
$y_0=0$.

{\it Proof of Theorem 1}. (a) $\Longrightarrow$ (b).
Suppose the condition (b) is violated and $n\in \{0,\dots,T-1\}$
is the largest number such that $W_n=\emptyset$ on the set
$$ A_n=\{\omega: \ri G_n^*\cap\ri Y_n=\emptyset\}=\{\omega:
0\not\in\ri (G_n^*-Y_n)\}$$ of positive measure. Note that the cone
$G_n^*-Y_n$ is not a subspace on this set. Therefore,
the mapping
$$D_n=(G_n^*-Y_n)^*\backslash (G_n^*-Y_n)^\perp I_{A_n}+
      (G_n^*-Y_n)^* I_{\Omega\backslash A_n}$$
take nonempty values.

According to standard measurable selection results (see e.g.
\cite[Theorem 5.2]{H75}) there exists an element
$$ -x_n\in L^0(D_n,\mathcal F_n).$$
In addition, $x_n\in L^0(-G_n,\mathcal F_n)$ and
$$ x_n\in L^0(Y_n^*,\mathcal F_n)=
  L^0(\mathcal K(W_{n+1},\mathcal F_n)^*,\mathcal F_n)=
  L^0(W_{n+1}^*,\mathcal F_n) $$
by Lemma 1. Now Lemma 3 implies that
$x_n\in\sum_{t=n+1}^T L^0(G_t,\mathcal F_t)$ and there exists
a collection of elements $x_t\in L^0(-G_t,\mathcal F_t)$, $t=n,\dots,T$,
satisfying the condition $\sum_{t=n}^T x_t=0$.

Let us show that this collection, augmented by the elements
$x_t=0$, $t=0,\dots,n-1$, violates the condition (a).
Suppose $x_t\in G_t\cap(-G_t)$ a.s., $t=0,\dots,T.$ By Lemma 5
we have
$$x_n=-\sum_{t=n+1}^T x_t\in
  L^0(Y_n^\perp,\mathcal F_n),$$
and, consequently, $-x_n\in (G_n^*-Y_n)^\perp$ a.s.,
in contradiction to the definition of $x_n$.

(b) $\Longrightarrow$ (a). Suppose $x_t\in L^0(-G_t,\mathcal F_t)$,
$\sum_{t=0}^T x_t=0$ is a collection of elements, violating condition (a).
Assuming also that condition (b) is satisfied, we shall
arrive at a contradiction.

Let $x_n\not\in L^0(G_n\cap (-G_n),\mathcal F_n)$.
Then $x_n\not\in (G_n^*)^\perp$ on a set of positive measure
$A_n\in\mathcal F_n$, while
$$ -x_n\in G_n=(G_n^*)^*\subset W_n^*\ \  {\rm a.s.}$$

If $x_n\in W_n^\perp$ on a subset $A_n'\subset A_n$ of positive
measure, then $x_n$ properly separates the sets $W_n$ and $G_n^*$
for $\omega\in A_n'$. By Lemma 6 this is impossible, since
$W_n\neq\emptyset$ a.s.

Now assume that $-x_n\in W_n^*\backslash W_n^\perp$ a.s. on $A_n$.
Then
\begin{equation} \label{1.6}
y_n=-\sum_{t=n}^T x_t\not\in L^0(W_n^\perp,\mathcal F_{n-1}).
\end{equation}
For $n=T$ this is evident, while for $n<T$ this is implied by
the equality $y_n=-x_n+y_{n+1}$ and Lemma 4, stating that
$$y_{n+1}=-\sum_{t=n+1}^T x_t\in L^0(Y_n^*,\mathcal F_n)\subset
  L^0(W_n^*,\mathcal F_n).$$

We get the contradiction, since by Lemma 7 the relation
(\ref{1.6}) is false under condition (b). The proof is comlpete.

{\bf 2. No-arbitrage criteria for the market with transaction costs.}
\setcounter{section}{2}
\setcounter{equation}{0}
In the sequel we consider the model of the market with transaction costs,
proposed in the paper \cite{S04}, where its economical meaning
and the connection with the models of \cite{K99}, \cite{KS01}, \cite{KRS02}
are discussed.

Assume that there are $d$ assets (currencies, as an example). Denote
by $\pi^{ij}$ the number of units of asset $i$, which can be
exchanged for one unit of asset $j$. Then one unit of asset $j$
can be exchanged for $1/\pi^{ji}$ units of asset $i$.
In other words, $1/\pi^{ji}$, $\pi^{ij}$ are the
bid and ask prices of asset $j$ in terms of asset $i$.

The square matrix $\Pi=(\pi^{ij})_{i,j=1}^d$ is called a
{\it bid-ask matrix} if
\begin{itemize}
\item[(i)] $\pi^{ij}>0,\ 1\le i,j\le d$;
\item[(ii)] $\pi^{ii}=1,\ 1\le i\le d$;
\item[(iii)] $\pi^{ij}\le\pi^{ik}\pi^{kj},\ 1\le i,j,k\le d$.
\end{itemize}
Condition (iii) means that the direct exchange is not worse than
any chain of exchanges. The interpretation of
conditions (i), (ii) is obvious.

Portfolios are the vectors in $\mathbf R^d$,
whose components describe the number of units of
each asset, hold by an investor. The convex cone $K(\Pi)$, spanned
by the vectors $\{e_i\}_{i=1}^d$
of the standard basis in $\mathbf R^d$, and the vectors
$\pi^{ij} e_i-e_j$ is called a {\it solvency cone}.
The conjugate cone and its relative interior
have the form (see \cite{S04}):
$$K^*(\Pi)=\{w\in\mathbf R^d_+:\pi^{ij}w_i-w_j\ge 0,
  \ 1\le i,j\le d\},$$
$$\ri K^*(\Pi)=\left\{w\in\intern\mathbf R^d_+:\frac{w_j}{w_i}\in
  \ri\left[\frac{1}{\pi^{ji}},\pi^{ij}\right],
  \ 1\le i,j\le d\right\},$$
where $\mathbf R^d_+=\{w\in\mathbf R^d:w_i\ge 0,\ i=1,\dots,d\}$.

The cone $K(\Pi)$ defines the preference relation
on the set of all portfolios, since from investor's point of view
elements of this cone are not worse than the zero portfolio.
Actually, every portfolio $\{\pi^{ij} e_i-e_j\}$
consists of $\pi^{ij}$ units long position of asset $i$ and
one unit short position of asset $j$.
Clearly, it can be transformed to the zero portfolio
(can be liquidated).

Assume that the asset price evolution is described
by an adapted to $(\mathcal F_t)_{t=0}^T$ stochastic process
$(\Pi_t)_{t=0}^T$, taking values in the set of
bid-ask matrices. Denote by $K_t=K(\Pi_t)$
the correspondent process of solvency cones.

Let an $\mathcal F_t$-measurable $d$-dimensional random
vector $\theta_t$ describe investor's portfolio at time $t$.
The portfolio process (shortly, portfolio)
$(\theta_t)_{t=0}^T$ is called a {\it self-financing} if
$$ \theta_t-\theta_{t-1}\in L^0(-K_t,\mathcal F_t),\ \ t=0,\dots,T,$$
where $\theta_{-1}=0$. At that, the portfolio $\theta_t$ at time $t$
is dominated by the portfolio $\theta_{t-1}$, hold in the previous
time moment.

Denote by $A_t(\Pi)$ the convex cone in $L^0(\mathbf R^d,\mathcal F_t)$,
formed by the elements $\theta_t$ of all self-financing portfolios
$\theta$:
$$ A_t(\Pi)=\sum_{j=0}^t L^0(-K_j,\mathcal F_j).$$

By the definition of \cite{S04}, the bid-ask process $(\Pi_t)_{t=0}^T$
satisfies the {\it robust no-arbitrage} (NA$^r$) condition
if there exists a bid-ask process $(\widetilde\Pi_t)_{t=0}^T$ such that
\begin{equation} \label{2.1}
[1/\widetilde\pi_t^{ji},\widetilde\pi_t^{ij}]\subset\ri
[1/\pi_t^{ji},\pi_t^{ij}]
\end{equation}
for all $i,j, t$ and
$$ A_T(\widetilde\Pi)\cap L^0(\mathbf R^d_+,\mathcal F_T)=\{0\}.$$

This condition requires the existence of a process $\widetilde\Pi$
with smaller, compared to $\Pi$, bid-ask spreads, and such that
it is arbitrage-free in a traditional un\-der\-stan\-ding.
The last statement means that it is impossible to "get
a riskless profit", i.e. the absence of a self-financing
portfolio (arbitrage strategy), which at time $T$ have
non-negative components and is non-trivial.
An additional stipulation is that the zero bid-ask spread is
regarded as "smaller than itself".

\begin{teo} For a bid-ask process $(\Pi_t)_{t=0}^T$
the following conditions are equivalent:
\begin{itemize}
\item[\rm (a)] NA$^r$ condition is satisfied;
\item[\rm (b)] there exists a $d$-dimensional $\mathbf P$-martingale
           $Z=(Z_t)_{t=0}^T$ such that
           $Z_t\in L^0(\ri K_t^*,\mathcal F_t)$, $t=0,\dots,T$;
\item[\rm (c)] all elements $W_t$ of the sequence of the set-valued maps
 $$W_T=K_T^*,$$
 $$W_t=\cl(\ri K_t^*\cap\ri Y_t),\ \ \
   Y_t=\cl(\conv\mathcal K(W_{t+1},\mathcal F_t)), \ \ 0\le t\le T-1$$
take nonempty values a.s.;
\item[\rm (d)] if $\sum_{t=0}^T x_t=0$, where
$x_t\in L^0(-K_t,\mathcal F_t)$, then $x_t=L^0(K_t\cap(-K_t),\mathcal F_t)$,
$t=0,\dots,T$.
\end{itemize}
\end{teo}

Here condition (c) is new. Its equivalence to condition (d)
is established in Theorem 1. Note, however, that the assertion of
Theorem 1 is more general. In particular, the cones $G_t(\omega)$
are not assumed to be polyhedral in this theorem.

The equivalence of conditions (a) and (b) was established in
the paper \cite{S04} (Theorem 1.7). The implications
(a) $\Longrightarrow$ (d) and (d) $\Longrightarrow$ (b) were proved
in \cite{KRS03} (see Lemma 5, Corollary 1, and also
the comment after it).

We emphasize, that for the cones $K_t^*$ with nonempty interior
(i.e. under the {\it efficient friction} condition \cite{KRS02}),
the equivalence of (b) and (c) is implied by the results of
\cite{R05}.

Note also that, putting
$$\theta_{-1}=0,\ \ \theta_t=\theta_{t-1}+x_t,\ \
0\le t\le T,$$
condition (d) can be reformulated as follows:
for any self-financing portfolio $\theta$
the equality $\theta_T=0$ implies that
$$\theta_t-\theta_{t-1}\in L^0(K_t\cap (-K_t),\mathcal F_t),
\ \ 0\le t\le T.$$

{\bf 3. Market model with a bank account.}
\setcounter{section}{3}
\setcounter{equation}{0}
There is no preffered asset in the model, considered in Section 2.
This is quite suitable for the description of a currency market.
At the same time, the asset exchanges are often cannot be done directly,
but only by buying and selling of assets with the use of a "canonical"
currency. In the present section we consider this situation as
a partial case of the general model.

So, consider the market model with a preffered asset (bank account).
Further, assume that the exchange operation of asset $i$
for asset $j$ is reduced to the exchange of asset $i$ for
the bank account and the exchange of the obtained money for
asset $j$. If the bank account is the first asset, then
the elements of the bid-ask matrix are represented in the form
$\pi^{ij}=\pi^{i1}\pi^{1j}$ \cite{G05}. At that
$$S^{b,i}=1/\pi^{i,1},\ \ S^{a,i}=\pi^{1,i},\ \ i=2,\dots,d$$
are the bid and ask prices of asset $i$, expressed in the units
of the bank account.

Conversely, if the $d$-dimensional vectors $S^b$, $S^a$ are such that
$0<S^{b,i}\le S^{a,i}$ and $S^{b,1}=S^{a,1}=1$; then the matrix
$\pi^{ij}=S^{a,j}/S^{b,i}$, $i\neq j$; $\pi^{ii}=1$
satisfies the conditions (i) -- (iii).
Let us call it a bid-ask matrix {\it for a market model with a bank
account}. It is easy to check that in such a model the solvency
cone $K(\Pi)$ is spanned by the vectors
$\{e_i-S^{b,i}e_1\}_{i=2}^d$, $\{S^{a,j}e_1-e_j\}_{j=2}^d$,
$\{e_i\}_{i=1}^d$. The conjugate cone and its relative interior
are described as follows
\begin{equation} \label{3.1}
K^*(\Pi)=\{w\in\mathbf R^d_+: S^{b,i} w_1\le w_i\le S^{a,i}w_1,\
i=2,\dots,d\},
\end{equation}
\begin{equation} \label{3.2}
\ri K^*(\Pi)=\{w\in\intern\mathbf R^d_+:
w_i/w_1\in\ri [S^{b,i},S^{a,i}],\ i=2,\dots,d\}.
\end{equation}

To establish the connection between the set
$$ C=[S^{b,2},S^{a,2}]\times\dots\times[S^{b,d},S^{a,d}]$$
and these cones it is convenient to utilize the representation
of points $x$ in $\mathbf R^{d-1}$ as the rays $\{\lambda (1,x):
\lambda\ge 0\}\subset \mathbf R_+\times
\mathbf R^{d-1}$ (cf. \cite[p.~76]{R73}).
To this end we introduce the set-valued map
$$ x\mapsto H(x)=\{\lambda (1,x): \lambda\ge 0\},
$$
which sometimes is called the {\it H\"ormander transform} \cite{KK02}.
For any set $A\subset\mathbf R^{d-1}$ we have
$$  H(A)=\bigcup_{x\in A} H(x)=
   \{\lambda (1,x): \lambda\ge 0,\ x\in A\}=
   \{(\mu,y):\mu\ge 0,\ y\in\mu A\}.$$

Let us also introduce the set-valued map  $x\mapsto\widehat{ H}(x)=
\{\lambda (1,x): \lambda>0\}$. Evidently,
\begin{equation} \label{3.3}
  H\left(\bigcup_{i=1}^\infty A_i\right)=
  \bigcup_{i=1}^\infty H(A_i),\ \ \
  H(\conv A)=\conv H(A),
\end{equation}
\begin{equation} \label{3.4}
 \widehat H(A\cap B)=\widehat H(A)
 \cap\,\widehat H(B).
\end{equation}
Moreover, $  H(\cl A)\subset \cl  H(A).$ Thus
\begin{equation} \label{3.5}
\cl H(\cl A)=\cl H(A).
\end{equation}
If, in addition, the set $A$ is bounded then, clearly,
the set $ H(\cl A)$ is closed and
\begin{equation} \label{3.6}
  H(\cl A)=\cl H(A)=
 \cl\widehat H(A).
\end{equation}
Note also that if the set  $A$ is convex then, by Corollary 6.8.1 of
\cite{R73}, we have
\begin{equation} \label{3.7}
 \ri\widehat H(A)=\ri H(A)= \widehat H(\ri A).
\end{equation}

The following result shows that for a set valued map $F$ the
closure of the H\"ormander transform $\cl H$ commutes with the
procedure of calculation of the regular conditional upper
distribution.
\begin{lem} Let $F$ be a closed $\mathcal F$-measurable
set-valued map and $\mathcal H\subset\mathcal F$.
Then the mapping $\omega\mapsto\cl H(F(\omega))$
is $\mathcal F$-measurable and
\begin{equation} \label{3.8}
\mathcal K(\cl H(F),\mathcal H)=
\cl H (\mathcal K(F,\mathcal H)).
\end{equation}
\end{lem}
{\it Proof.} Consider a Castaing representation $\{\xi_i\}_{i=1}^\infty$
for $F$ and let $\{\lambda_j\}_{j=1}^\infty$ be a countable dense
subset of $\mathbf R_+$. Clearly,
$\{(\lambda_j,\lambda_j\xi_i)\}_{i,j=1}^\infty$ is the Castaing
representation for the mapping $\cl H(F)$. Hence, this mapping is
$\mathcal F$-measurable.

It is sufficient to consider the case $F(\omega)\neq
\emptyset$ a.s., since otherwise the sets in the both sides of
equality (\ref{3.8}) are empty by the definition of the
support of the regular conditional upper distribution.
Using formula (\ref{1.2}), we get
\begin{eqnarray*}
\mathcal K(\cl H(F),\mathcal H) &=&
\cl\left(\bigcup_{i,j=1}^\infty\mathcal K((\lambda_j,\lambda_j\xi_i),
\mathcal H)\right)=\cl\left(\bigcup_{i,j=1}^\infty\lambda_j
\mathcal K((1,\xi_i),\mathcal H)\right)\\
&=&\cl\left(\bigcup_{i=1}^\infty\bigcup_{j=1}^\infty\lambda_j
(1,\mathcal K(\xi_i,\mathcal H))\right)=
\cl\left(\bigcup_{i=1}^\infty H(\mathcal K(\xi_i,\mathcal H))\right).
\end{eqnarray*}
In the last equality we take into account that the set
$\bigcup_{j=1}^\infty\lambda_j (1,\mathcal K(\xi_i,\mathcal H))$
is dense in $H(\mathcal K(\xi_i,\mathcal H))$ for any $i$.

On the other hand, formulas (\ref{1.2}), (\ref{3.5}), (\ref{3.3})
imply that
\begin{eqnarray*}
  \cl H (\mathcal K(F,\mathcal H))&=&
  \cl H\left(\cl\left(\bigcup_{i=1}^\infty
  \mathcal K(\xi_i,\mathcal H)\right)\right)=
  \cl H\left(\bigcup_{i=1}^\infty
  \mathcal K(\xi_i,\mathcal H)\right)\\
  &=&
  \cl\left(\bigcup_{i=1}^\infty
  H(\mathcal K(\xi_i,\mathcal H))\right).
\end{eqnarray*}
This completes the proof.

Recurring to the market model with a bank account, consider
an adapted $d$-dimensional stochastic processes
$S^b=(S_t^b)_{t=0}^T$, $S^a=(S_t^a)_{t=0}^T$,
satisfying the conditions
$$0<S_t^{b,i}\le S_t^{a,i},\ \ i=2,\dots,d;\ \ \
  S^{b,1}_t=S^{a,1}_t=1$$
for all $t$. Put
$$C_t=[S_t^{b,2},S_t^{a,2}]\times\dots\times[S_t^{b,d},S_t^{a,d}].$$

Let $\widehat w=(w_2,\dots,w_d)$. Then
\begin{equation} \label{3.9}
 H(C_t)=\{(w_1,\widehat w): w_1\ge 0,\
 \widehat w\in w_1 C_t\}=K^*_t,\ \ K_t=K(\Pi_t);
\end{equation}
\begin{equation} \label{3.10}
 \widehat H(\ri C_t)=\{(w_1,\widehat w): w_1> 0,\
 \widehat w\in w_1 \ri C_t\}=\ri K^*_t
\end{equation}
by (\ref{3.1}), (\ref{3.2}).

\begin{teo}
Let  $(\Pi_t)_{t=0}^T$ be a bid-ask process
for the model with a bank account. Then the following
conditions are equivalent:
\begin{itemize}
\item[\rm (a)] NA$^r$ condition is satisfied;
\item[\rm (b)] there exist an equivalent to $\mathbf P$ probability measure
              $\mathbf Q$ and a $d-1$-dimensional $\mathbf Q$-martingale
              $S=(S_t)_{t=0}^T$, $S_t=(S_t^2,\dots,S_T^d)$ such that
              $$S_t^i\in L^0(\ri [S_t^{b,i},S_t^{a,i}],\mathcal F_t),
              \ \ i=2,\dots,d,\ t=0,\dots,T;$$
\item[\rm (c)] all elements $V_t$ of the sequence of set-valued maps
 $$V_T=C_T,$$
 $$V_t=\cl(\ri C_t\cap\ri X_t),\ \ \
   X_t=\cl(\conv\mathcal K(V_{t+1},\mathcal F_t)), \ \ 0\le t\le T-1$$
take nonempty values a.s.
\end{itemize}
\end{teo}

{\it Proof.} Although only condition (c) is new, for the sake
of completeness, we prove that conditions (b) and (c) of Theorem 3
are equivalent to the correspondent conditions of Theorem 2.

Suppose, condition (b) of Theorem 3 is satisfied. Following \cite{S04},
we introduce the stochastic process $Z=(Z_t)_{t=0}^T$ by
the formulas
$$ Z_t^1=\mathbf E\left(\frac{d\mathbf Q}{d\mathbf P}\biggr|
   \mathcal F_t\right),
   \ \ \ Z_t^i=Z_t^1 S_t^i.$$
This process is a $\mathbf P$-martingale and
$$ Z_t=(Z_t^1,Z_t^1 S_t^i)\in\widehat H(\ri C_t)=\ri K_t^*\ {\rm a.s.}$$
Consequently, condition (b) of Theorem 2 holds true.

Now suppose that condition (b) of Theorem 2 is satisfied and
$Z$ is the process, mentioned in this condition. Again, following
\cite{S04}, we introduce an equivalent to $\mathbf P$ probability measure
$\mathbf Q$ with the density $d\mathbf Q/d\mathbf P=Z_T^1/Z_0^1$ and
$\mathbf Q$-martingale $S_t=(Z_t^2/Z_t^1,\dots,Z_t^d/Z_t^1)$.
Evidently, $S_t\in \ri C_t$  a.s. for all $t$ and
condition (b) of Theorem 3 holds true.

The equivalence of conditions (c) of Theorems 2 and 3 is implied
by the equalities
\begin{equation} \label{3.11}
W_t=H(V_t),\ \ \ 0\le t\le T,
\end{equation}
which are established below.

According to formula (\ref{3.9}) we have
$$ H (V_T)=H (C_T)=K_T^*=W_T.$$
Let $W_{t+1}=H(V_{t+1})$. Keeping in mind that the sets
$C_t(\omega)$ are bounded, we use (\ref{3.6}) and also
formulas (\ref{3.4}), (\ref{3.10}), (\ref{3.7}):
\begin{eqnarray*}
H(V_t) &=& \cl H (\ri C_t\cap\ri X_t)
=\cl\widehat H(\ri C_t\cap\ri X_t)=
 \cl(\widehat H(\ri C_t)\cap \widehat H(\ri X_t)) \\
&=& \cl (\ri K_t^*\cap \ri H(X_t)).
\end{eqnarray*}
Furthermore, by the definition of $X_t$ and formulas (\ref{3.5}), (\ref{3.3}),
(\ref{3.8}) (taking into consideration the boundedness of $V_{t+1}(\omega)$
and formula (\ref{3.6})), we get
\begin{eqnarray} \label{3.12}
\cl H (X_t)&=&\cl H (\cl(\conv\mathcal K(V_{t+1},\mathcal F_t)))=
\cl H (\conv\mathcal K(V_{t+1},\mathcal F_t))\nonumber\\
&=& \cl(\conv H(\mathcal K(V_{t+1},\mathcal F_t)))
=\cl(\conv(\cl H(\mathcal K(V_{t+1},\mathcal F_t))))\nonumber\\
&=&\cl(\conv\mathcal K(H(V_{t+1}),\mathcal F_t))=
\cl(\conv\mathcal K(W_{t+1},\mathcal F_t))\nonumber\\
&=& Y_t.
\end{eqnarray}
Thus, $H(V_t)=\cl (\ri K_t^*\cap \ri Y_t)=W_t$ and equalities
(\ref{3.11}) hold true. The proof is complete.

{\bf 4. On construction of arbitrage strategies.}
Assume that condition (c) of Theorem 2 is violated and
$n$ is the largest number such that $W_n=\emptyset$ on a set $A_n$
of positive measure. Following the argumentation, used in the proof
of the implication (a) $\Longrightarrow$ (b) of Theorem 1 (for $G_t=K_t$),
consider an arbitrary element
$$ -x_n\in L^0((K_n^*-Y_n)^*\backslash (K_n^*-Y_n)^\perp I_{A_n}+
      (K_n^*-Y_n)^* I_{\Omega\backslash A_n},\mathcal F_n).$$
This element is characterized by the property
$x_n\in L^0(-K_n\cap Y_n^*,\mathcal F_n)$ and by the fact that
$x_n(\omega)$ properly separates the sets $K_n^*(\omega)$ and
$Y_n(\omega)$ for almost all $\omega\in A_n$.

The element $x_n$ admits a representation
$$x_n=-\sum_{t={n+1}}^T x_t,\ \ \ x_t\in L^0(-K_t,\mathcal F_t)$$
and the collection $\{x_t\}_{t=0}^T$, where $x_t=0$, $t<n$,
violates condition (d) of Theorem 2 (see the proof of Theorem 1).
Hence, there exists a number $m\ge n$ such that $x_m\not\in
F_m=K_m\cap (-K_m)$ on a set $B_m\in\mathcal F_m$ of positive
measure.

Let $\widetilde\Pi$ be a bid-ask process, satisying (\ref{2.1}).
Then $\intern (-K(\widetilde\Pi_m))\supset -K_m\backslash F_m$
\cite[Lemma 2.5]{S04} and there exists a random vector $\varepsilon_m\in
L^0(\mathbf R^d_+,\mathcal F_m)$ not equal to zero on $B_m$ and such that
$x_m+\varepsilon_m\in -K(\widetilde\Pi_m)$ a.s.

Put $\varepsilon_t=0$, $t\neq m$ and introduce the portfolio process
$\theta$ by the formulas
$$\theta_{-1}=0,\ \ \theta_t=\theta_{t-1}+x_t+\varepsilon_t,\ \
0\le t\le T.$$
Clearly, $\theta$ is an arbitrage self-financing portfolio for the process
$\widetilde{\Pi}$:
$$\theta_T=\sum_{t=0}^T (\theta_t-\theta_{t-1})=\varepsilon_m\in
A_T(\widetilde\Pi)\cap (L^0(\mathbf R^d_+,\mathcal F_T)\backslash\{0\}).$$

{\bf Example.} Suppose $\Omega=\{\omega_1,\omega_2\}$, the algebra
$\mathcal F$ is generated by one-point sets,
$\mathbf P(\omega_i)>0$, $i=1,2$ and $\mathcal F_0=\mathcal F_1=
\{\emptyset,\Omega\}$, $\mathcal F_2=\mathcal F$.
Consider the model with a bank account (for $T=2$, $d=3$),
where the sets $C_t$ are of the following form:
$$ C_0=[2,6]\times[2,6],\ \ \ C_1=[4,8]\times[4,8];$$
$$ C_2(\omega_1)=[9,9]\times[6,6],\ \ \
   C_2(\omega_2)=[4,4]\times[1,1].$$

By formula (\ref{1.1}) we have
$$\mathcal K(V_2,\mathcal F_1)=\mathcal K(C_2,\mathcal F_1)=
C_2(\omega_1)\cup C_2(\omega_2).$$
Consequently, the set $X_1=\cl(\conv \mathcal K(V_2,\mathcal F_1))$
is the line segment, connecting the points $(4,1)$ and $(9,6)$;
the set $V_1=\cl(\ri C_1\cap\ri X_1)$ is the line segment,
connecting the points $(7,4)$ and $(8,5)$ (see Fig.\,1).
Further, we get
$$X_0=\cl(\conv \mathcal K(V_1,\mathcal F_0))=V_1, \ \
  V_0=\cl(\ri C_0\cap\ri X_0)=\emptyset.$$

By Theorem 3 the NA$^r$ condition is violated and there exists
an arbitrage strategy on the market under consideration.
To construct such a strategy, according to the scheme
stated above, at first we ought to find an element $x_0\in -K_0$,
properly separating the cones $K_0^*=H(C_0)$ and
$Y_0=\cl H(X_0)=H(X_0)$ (here formulas (\ref{3.9}),
(\ref{3.12}) are used).

Note that if an element $p\in\mathbf R^{d-1}$ properly separates
the sets $A,B\subset\mathbf R^{d-1}$ and
$$ \langle p,\xi\rangle\le
   \alpha\le \langle p,\eta\rangle,\ \ \ \xi\in A,\ \eta\in B,$$
then the element $-x=(\alpha,-p)\in H(A)^*$
properly separates the cones $H(A)$, $H(B)$.

Since $p=(1,0)$ properly separates the sets $C_0$, $X_0$
and $\langle p,\eta\rangle\ge 7$, $\eta\in X_0$, we conclude that
the element $x_0=(-7,1,0)$ is desirable.

The next step consists in the obtaining of a representation $-x_0=x_1+x_2$,
$x_1\in -K_1$, $x_2\in -K_2(\omega)$ (the vector $x_2$ is
$\mathcal F_1$-measurable, i.e. constant). Let's find the
conjugate cones:
$$ K_1^*=H(C_1)=
   \cone\{(1,4,4),(1,8,4),(1,4,8),(1,8,8)\};$$
$$ K_2^*(\omega_1)=H(C_2(\omega_1))=\cone\{(1,9,6)\},$$
$$ K_2^*(\omega_2)=H(C_2(\omega_2))=\cone\{(1,4,1)\}.$$
We get the following system of linear equalities and inequalities
for $x_i=(x_i^1,x_i^2,x_i^3)$, $i=1,2$:
$$ x_1^1+x_2^1=7,\ \ x_1^2+x_2^2=-1,\ \ x_1^3+x_2^3=0$$
$$ x_1^1+4 x_1^2+4 x_1^3\le 0, \ \ x_1^1+8 x_1^2+4 x_1^3\le 0$$
$$ x_1^1+4 x_1^2+8 x_1^3\le 0, \ \ x_1^1+8 x_1^2+8 x_1^3\le 0$$
$$ x_2^1+9 x_2^2+6 x_2^3\le 0, \ \ x_2^1+4 x_2^2+x_2^3\le 0$$
It is not difficult to check that
$x_1=(4,0,-1)$, $x_2=(3,-1,1)$
is an unique solution of this system.

Assume that the prices of assets are expressed in conventional units
(c.u.). The strategy $\theta$: $\theta_{-1}=0$,
$\theta_t-\theta_{t-1}=x_t$, $t=0,1,2$ has the following
interpretation. For $t=0$ the investor buys one unit
of the second asset at the price 6 c.u. and "throws away"
1 c.u. For $t=1$ he sells short one unit of the third asset,
and gets 4 c.u. on his bank account. For $t=2$ he sells
one unit of the second asset and buys one unit of the third one. At
that, he gets more 3 c.u. on his bank account. As a result,
his portfolio $\theta_2$ is zero.

Since for $t=0$ it is not necessary to throw away 1 c.u.,
this portfolio can easily be transformed to an arbitrage one,
by changing  $x_0$ to $\widetilde{x}_0=(-6,1,0)\in-K_0$,
$$ \widetilde{x}_0+x_1+x_2
=(1,0,0) \in A_2(\Pi)\cap(\mathbf R^3\backslash\{0\}).$$
All the more, this portfolio still produces arbitrage for a smaller
bid-ask spread.

We also mention that the obtained arbitrage strategy involve
three steps: all increments $x_t$, $t=0,1,2$ of the portfolio
are non-zero. There are no two-step arbitrage strategies in the examined
example. Indeed, Fig.\,1 shows that when one of the steps is omitted,
condition (c) of Theorem 3 holds true.

\newpage
\begin{picture}(12,12)
\put(0,0){\vector(1,0){10}}
\put(0,0){\vector(0,1){10}}
\put(2,2){\line(1,0){4}}
\put(2,2){\line(0,1){4}}
\put(6,2){\line(0,1){4}}
\put(2,6){\line(1,0){4}}
\put(4,4){\line(1,0){4}}
\put(4,4){\line(0,1){4}}
\put(8,4){\line(0,1){4}}
\put(4,8){\line(1,0){4}}
\multiput(4,1)(0.5,0.5){9}{\line(1,1){0.4}}
\put(8.55,5.55){\line(1,1){0.45}}
\put(9.2,6){$C_2(\omega_1)$}
\put(3.8,0.5){$C_2(\omega_2)$}
\put(2.2,2.2){$C_0$}
\put(7.4,7.5){$C_1$}
\put(8.35,2.8){\dashbox{0.2}(0.8,0.8){$V_1$}}
\put(7.6,4.4){\vector(1,-1){0.7}}
\put(6.9,3.5){$(7,4)$}
\put(8.2,4.8){$(8,5)$}
\put(7,4){\circle*{0.1}}
\put(8,5){\circle*{0.1}}
\put(9,6){\circle*{0.1}}
\put(4,1){\circle*{0.1}}
\put(-0.3,2){$2$}
\put(-0.3,4){$4$}
\put(-0.3,6){$6$}
\put(-0.3,8){$8$}
\put(2,-0.4){$2$}
\put(4,-0.4){$4$}
\put(6,-0.4){$6$}
\put(8,-0.4){$8$}
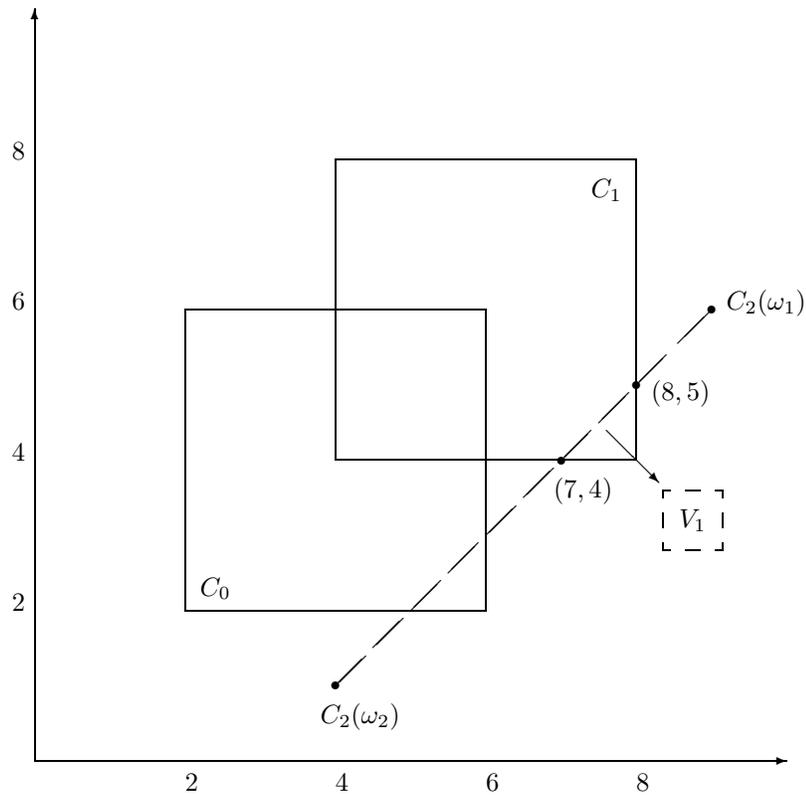
\put(0.5,-1.4){{\bf Fig.1.} Example of a market, admitting}
\put(1.7,-2){a three-step arbitrage strategy.}
\end{picture}
\end{document}